\documentclass[12pt]{article}

\usepackage{epsfig,psfrag,harvard,amssymb,fullpage,graphicx}

\citationmode{abbr}
\begin{document}
\def\sign{{\rm sign}}
\def\PP{{\cal P}}

\let\hat\widehat
\let\tilde\widetilde

\def\corr{{\rm corr}}
\def\min{{\rm min}}
\def\sN{ {\cal N}}
\def\sM{ {\cal M}}
\def\sA{{\cal A}}
\def\argmax{{\rm argmax}}
\title{Correlation-sharing for detection of  differential gene expression}

\author{ Robert  Tibshirani \thanks{Dept. of Health Research \& Policy, and 
Department of Statistics, Stanford
    University, Stanford, CA 94305. Email: {\tt
      tibs@stat.stanford.edu}.} \\ 
Larry Wasserman
  \thanks{Department of
    Statistics, Carnegie Mellon University, Pittsburg, PA. Email: {\tt
      larry@stat.cmu.edu}.}}
\maketitle
\vskip -.5in
\begin{abstract}
We propose a method for detecting differential gene expression that
exploits the  correlation between genes.
Our proposal averages the univariate scores of each feature with
the scores in correlation neighborhoods.   In a number of real and simulated
examples, the new method often exhibits lower false discovery rates
than  simple t-statistic
thresholding.
We also provide some  analysis of the asymptotic behavior of our proposal.
The general idea of correlation-sharing can be applied to other prediction
problems involving a large number of correlated features.
We give an example in protein mass spectrometry.
\end{abstract}
  
\section{Introduction}
We consider methods for detecting differentially expressed genes in
from a set of microarray experiments.
Consider the simple case of $m$ genes measured across two experimental
conditions.
A number of authors have proposed methods for detecting differential
gene expression,
including \citeasnoun{Ditetal00}, \citeasnoun{Newtonetal01} and \citeasnoun{Kerr2000}.
\citeasnoun{Storey2006}  presents an interesting, more general approach.

One widely used approach to this problem is
as follows. We compute a two-sample t-statistic $T_i$ for each gene, and then
call a gene significant if $|T_i|$ exceeds some threshold $c$.
Various values of $c$ are tried, using permutations of the sample labels 
to estimate the false discovery rate (FDR) for the procedure for each $c$.
A threshold $c$ is finally chosen based on the estimates of FDR and other
considerations, such as the ballpark number of significant genes that is
 desirable.
This  recipe roughly  describes the strategy used, for example,  in the Significance of
Microarrays (SAM) procedure
\cite{TTC01}.

In this paper we propose a simple method for potentially improving on the
thresholded t-statistic approach defined above.
The idea is to exploit correlation among the genes.
In a sense this general idea is not new, and exploratory
methods based on clustering have been proposed (e.g. \citeasnoun{Tibsetal2002}).
These methods require choices like the clustering metric and
linkage, and hence are somewhat subjective. 
The proposal presented here is much simpler, and hence it is easier to analyze
and assess its performance.

We start with t-statistics computed for each gene. Then we assign to each
gene a score $r_i$ equal to the average of all t-statistics for genes having correlation
at least $\rho(i)$ with that gene, choosing the best value of $\rho(i) \in [0,1]$
to maximize the average.
Finally, we 
call a gene significant if $|r_i|$ exceeds some threshold $c$.
The idea is that differentially expressed genes are likely to co-exist in
a pathway, and hence   will be correlated in our data. Hence use of the
score $r_i$ might provide a more accurate test of significance than 
that based on $t_i$. We call this approach ``correlation sharing"
Note that the choice $\rho(i)=1$ yields no sharing, giving $r_i=t_i$.
Hence the correlation-sharing method contains the thresholded t-statistic
approach as a special case.

As a motivating example, we generated data with 1000 genes and 30 samples.
The first 50 genes $i \in \PP=\{1,2,\ldots 50\}$ are generated as

\begin{equation}
  \label{eq:4}
  X_{ij}=Z_{ij}+ .75\cdot I(j>15)
\end{equation}
with  $Z_{ij} \sim N(0,1)$ and ${\rm corr}(Z_i,Z_{i'})+.0.8$,
where $Z_i=(Z_{i1},\ldots Z_{in})$
The remaining genes $x_{ij}, i>50$ were generated as $N(0,1)$.
The outcome variable  $Y_j$  equaled  2 for $16\leq j \leq 30$ and
1 otherwise.

\begin{figure}
\centerline{\epsfig{file=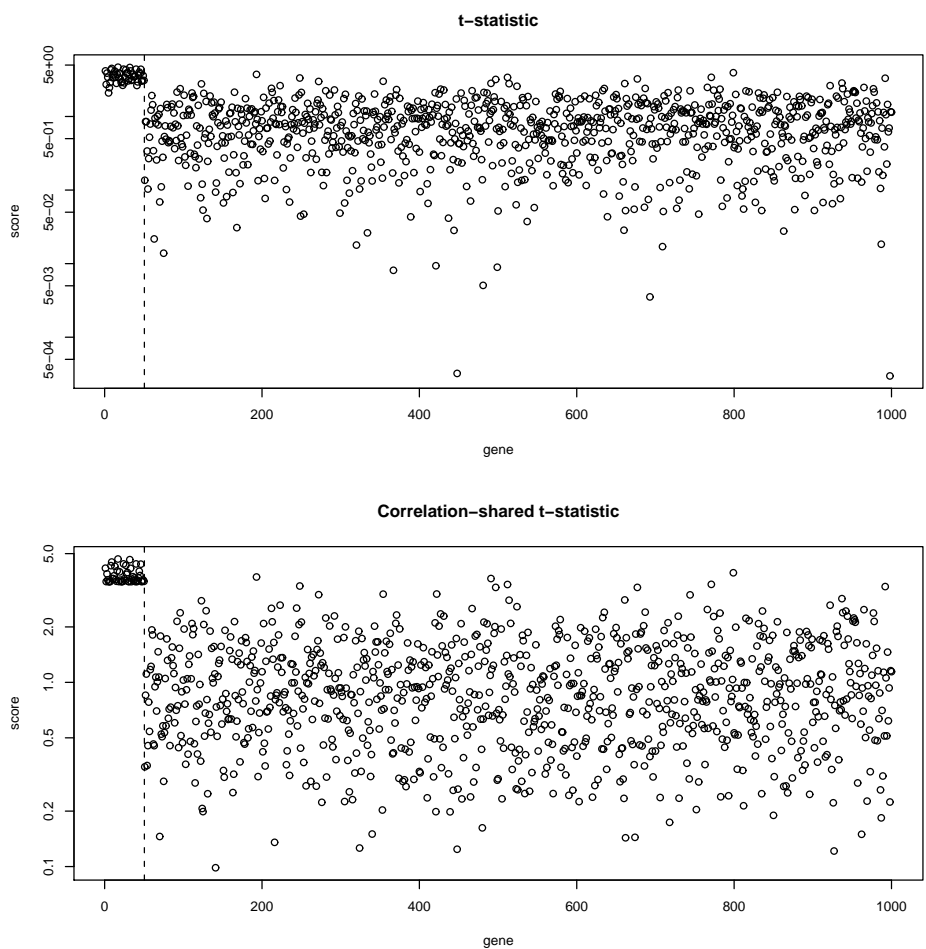,width=5in}}
\caption[curb0]{\small\em  T-statistics and correlation-shared
  T-statistics for simulated example 1.}
\label{curb0}
\end{figure}

Figure \ref{curb0} shows the
 t-statistics (top panel) and correlation-shared  t-statistics
(bottom panel). We see that in the bottom panel the scores
for the first 50 genes are  magnified.
This  leads to improved detection of the differentially expressed
genes, as we show in the next section.

The outline of this paper is a follows.
Section \ref{sec.def}  defines  correlation-sharing.
In section \ref{sec.res} we  discuss the concept of {\em residual correlation},
and its impact on correlation-sharing.
We apply our method to
four microarray  cancer datasets.
The skin data is examined more closely in section \ref{sec.skin}.
Some asymptotic results for correlation sharing are given
in section \ref{sec.theory}.
Section \ref{sec.protein} applies the method to a different kind of data--- 
protein mass spectra.
Finally in section \ref{sec.other} we discuss the application
of correlation sharing to other kinds of response variables,
and computational issues.

\section{Correlation sharing}
\label{sec.def}
Let $X$ be the $m\times n$ matrix of expression values, for $m$ genes
and $n$ samples. We assume that the samples fall into two groups $j=1$ and $2$.
We start with th standard (unpaired)  t-statistic
\begin{eqnarray}
 T_i= {{\bar x_{i2}- \bar x_{i1}}\over{s_i}}
\label{tstat}
\end{eqnarray}

Here $\bar x_{ij}$ is the   mean of gene $i$ in group $j$ and  $s_i=$ pooled within group standard deviation  of gene $i$.

Let $x_i$ denote the $ith$ row of $X$.
Define $C_{\rho}(i)=  \{ k: {\rm corr}(x_i, x_k) \geq \rho\}$, the    indices of the genes 
with correlation at least $\rho$  with gene $x_i$
Then we define 
\begin{eqnarray}
u_i&=&{\rm max}_{\{0 \leq \rho \leq 1\}}\;  {\rm ave}_{j \in C_{\rho}(i)} |T_j|\cr
r_i&=& {\rm sign}(T_i)\cdot u_i
\label{eq:corrshare}
\end{eqnarray}
We call this the ``correlation-shared'' t-statistic.
The method calls significant all genes having   $|r_i|>c$,  and estimates the false discovery rate (FDR) of the resultant gene list by permutations.
We vary  $c$  and examine the estimated FDR.

\begin{figure}
\begin{psfrags}
\psfrag{# genes called}{\footnotesize \hskip -.35in Number of genes called significant}
\psfrag{# genes called }{\footnotesize \hskip -.35in Number of genes called significant}
\psfrag{Number of false positive genes}{\footnotesize \hskip -.35in Number of false positive genes}
\psfrag{Number of false negative genes}{\footnotesize \hskip -.35in Number of false negative genes}
\centerline{\epsfig{file=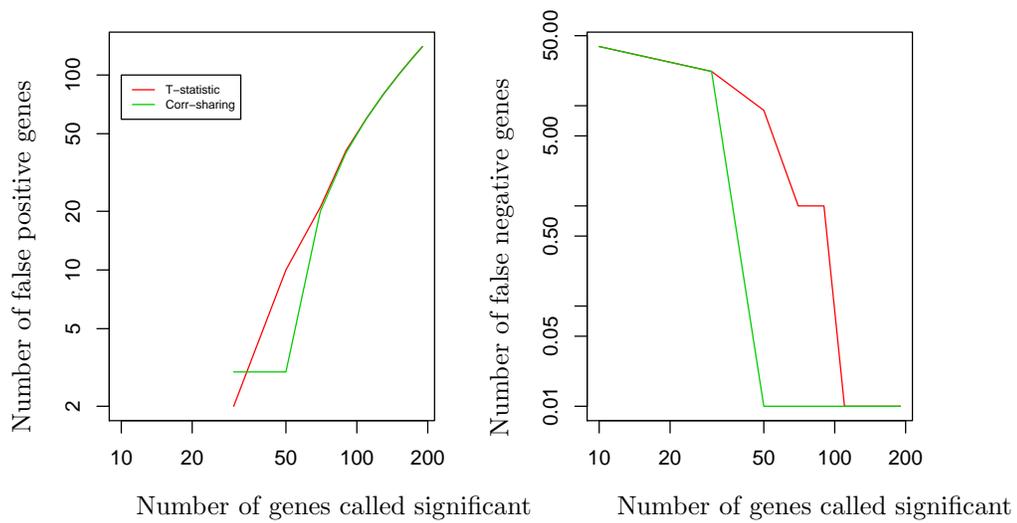,width=5in}}
\end{psfrags}
\caption[curb1]{\small\em Results for example 1. Left panel: Number of false positive genes versus number
of genes called significant. 
Right panel: Number of false negative genes versus number
of genes called significant.}
\label{curb1}
\end{figure}

Figure \ref{curb1} shows the results for  correlation
sharing
applied to the simulated data from model (\ref{eq:4}) . As the threshold is varied,
the number of genes called significant and the number of false positive
genes and false negative genes all change. 
We see that correlation sharing generally  yields  fewer false positive and
false negative genes
genes than  the t-statistic.

We  can also think of correlation-sharing as a method
for {\em supervised clustering}.
Let $\hat\rho(i)$ be the maximizing correlation for gene $i$,
from definition (\ref{eq:corrshare}).
Then the set of genes with indices  $C_{\hat p(i)}(i)$
is an adaptively chosen cluster, selected to maximize the average
``signal'' around gene $i$.
Unlike with most standard clustering methods, the 
 clusters  $C_{\hat p(i)}(i)$
are overlapping, rather than mutually disjoint.
We examine these clusters in  some examples later in this paper.

\begin{figure}
\begin{psfrags}
\psfrag{# genes called}{\footnotesize \hskip -.35in Number of genes called significant}
\psfrag{# genes called }{\footnotesize \hskip -.35in Number of genes called significant}
\psfrag{Number of false positive genes}{\footnotesize \hskip -.35in Number of false positive genes}
\psfrag{Number of false negative genes}{\footnotesize \hskip -.35in Number of false negative genes}
\centerline{\epsfig{file=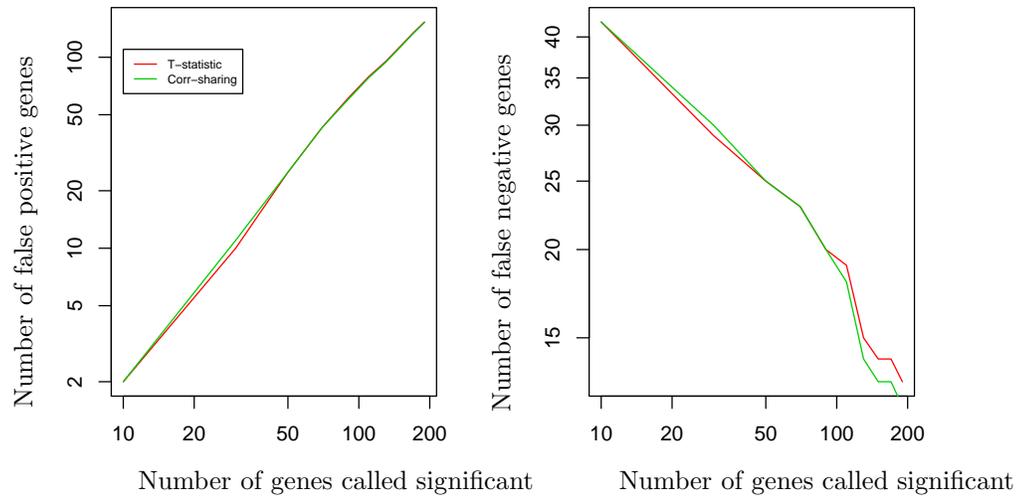,width=5in}}
\end{psfrags}
\caption[curb11]{\small\em Results for example 2. Here the non-null
genes have no  correlation before the group effect is added.}
\label{curb11}
\end{figure}

As a second example,  we  changed the data generation so that the first 50 genes had no correlation,
before the group effect was added.
Figure \ref{curb11} shows that the  advantage of correlation sharing has disappeared.

\section{Residual correlation among non-null and null genes}
\label{sec.res}
The previous example suggests that a key assumption in
for our proposal is that the  correlation between the  non-null genes
is  higher
than that for the null genes.

We need to say precisely what we mean by ``correlation''. Suppose for
a set of non-null genes $\PP$, the expression is $\beta$ units higher
in group $Y_j=2$ than it is in group $Y_j=1$:
\begin{eqnarray}
x_{ij} &=& \beta \cdot I(Y_j=2) +\varepsilon_{ij} \; {\rm for}\; i\in  \PP\cr
       &=&  \varepsilon_{ij} \; {\rm for} \;  i\notin  \PP 
\end{eqnarray}

Let $x_i= (x_{i1}, x_{i2}, \ldots x_{in})$. Then even if the errors
$\varepsilon_{ij}$ are all independent of one another, we have
$\corr(x_i, x_{i'})>0$ for $i, i'\in \PP$. That is, the treatment
effect induces an overall correlation between the genes in $\PP$.
However we would expect that the t-statistic would capture all of the
information
needed to decide if a gene is in $\PP$.

Instead, we  assume that there is {\em residual correlation} among
the genes in $\PP$:
\begin{eqnarray}
\corr(\varepsilon_i, \varepsilon_{i'}) >0; \; {\rm for} \; i,i' \in \PP
\end{eqnarray}
where $\varepsilon_i=(\varepsilon_{i1}, \ldots \varepsilon_{in})$.

For the simulated data of Figure \ref{curb0},
the estimated residual correlation is the correlation between
genes, after having removed the estimated effect of
treatment. Specifically, the  residual correlation is 
$\corr(x^*_{i}, x^*_{i'})$ where  $x^*_{ij}=x_{ij}-\hat x_{ij}$.
For the two sample case, for example,  $\hat x_{ij}=\bar x_{i2}-\bar x_{i1}$, $\bar
x_{ik}$ equaling
the average of $x_{ij}$ for samples in group $k$.

The average absolute  residual correlation for the non-null genes (the
first 50 genes)
equaled 0.47, while that for the null genes was 0.15, and the correlation
between the non-null and null genes was also 0.15.

Is there residual correlation in   real microarray data?  Biologically, genes will
be correlated if they are in the same pathway. However if that pathway
is not active in the  experimental conditions under study, the  genes in the
pathway will not
show large correlation. And the same genes will tend to be null, i.e.
will not differentially expressed in the experiment. The opposite
should be true
for differentially expressed genes.

To see if this assumption is reasonable in practice, we examine four 
microarray datasets: the skin data taken from \citeasnoun{rieger2004},
and Duke breast cancer data taken from \cite{huang2003},
the BRCA data taken from \citeasnoun{Hedetal2001} and the non-Hodgkins lymphoma
data from \citeasnoun{aR02}.
These are summarized in Table \ref{tab:summary}.

The false discovery rates of both the t-statistic and 
correlation-shared statistics depend on the total number of genes
input into the corresponding procedure.
Hence for fairness (and computational speed)
we started with the 2000 genes having largest overall variance
in each case.

To examine residual correlation,
we computed the two-sample t-statistics $T_i$ for each gene.
 Then we computed the average absolute  residual correlation for genes
satisfying $|T_i|>c$, with $c$ varying from the 99th to the 75 quantiles of
the  $|T_i|$ values. In the lymphoma data the outcome is survival time;
hence we instead computed the Cox's partial likelihood score statistic
for each gene (see section \ref{sec.other}).

\begin{figure}
\begin{psfrags}
\psfrag{Number of genes called}{\footnotesize \hskip -.35in Number of genes called significant}
\psfrag{ave abs corr}{\footnotesize \hskip -.35in Average absolute correlation}
\centerline{\epsfig{file=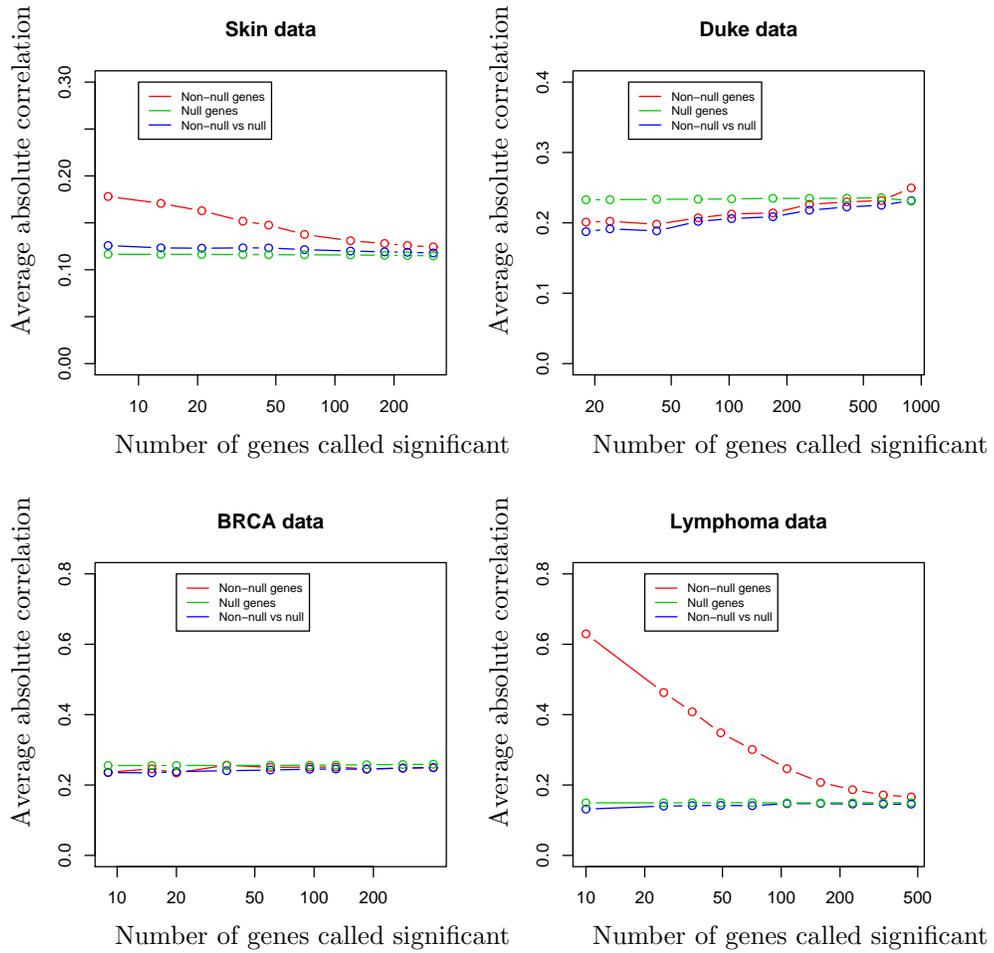,width=5in}}
\end{psfrags}
\caption[cancer.corr]{\small\em Average absolute residual correlation
  as a function of the number of genes called significant by the T or F-statistics.}
\label{cancer.corr}
\end{figure}

 The results are   shown in Figure \ref{cancer.corr}.
For the skin and lymphoma datasets data, the non-null genes have  higher correlation
with each other than they have with the null genes,
and also higher than that within the null-genes.
But for the Duke and BRCA2  datasets, this is not the case.

\begin{figure}
\begin{psfrags}
\psfrag{# called}{\footnotesize \hskip -.35in Number of genes called significant}
\psfrag{#FP}{\footnotesize \hskip -.35in Number of false positive genes}
\centerline{\epsfig{file=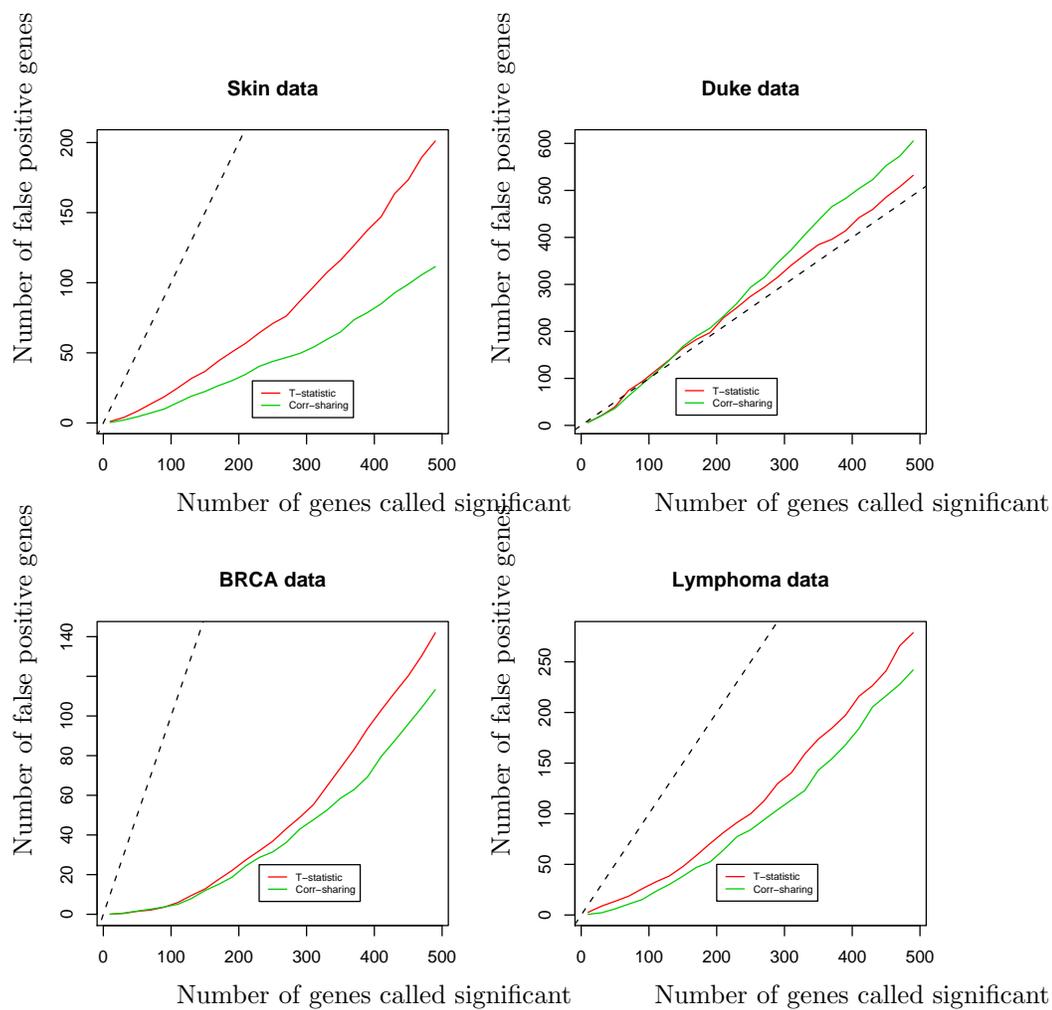,width=5in}}
\end{psfrags}
\caption[fpall]{\small\em
Results for four cancer datasets: plotted is the number of false positive genes
versus the number of genes called, for the standard t-statistic (red)
and the correlation-shared t-statistic (green).
The broken line is the $45^o$ line.}
\label{fpall}
\end{figure}

For the same four datasets,
 Figure \ref{fpall} shows
the estimated number of false positive genes is plotted against the
number of genes called significant, for both the t-statistic and  correlation shared t-statistic.
Correlation sharing exhibits lower FDR for all datasets except the Duke data,
where neither method does much as all.

\section{Skin data example}
\label{sec.skin}
We examine more closely the results for the skin data shown in the top left
panel of Figure \ref{fpall}.
 There are 12,625 genes and  58  patients:
44 normal patients and 14 with radiation sensitivity.

\begin{figure}
\centerline{\epsfig{file=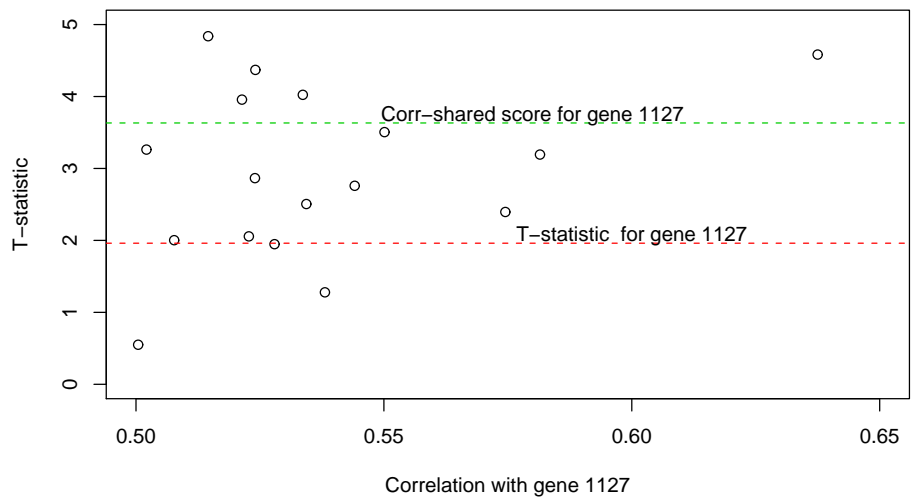,width=5in}}
\caption[skin3]{\small\em Skin data: a closer look at gene 1127}
\label{skin3}
\end{figure}

Figure \ref{skin3} illustrates how correlation sharing can magnify the
effect of a gene (\#1127 chosen as an example).
The figure shows all genes having correlation at least
0.5 with gene \# 1127.
Its raw t-statistic is about 2.0
Notice that the  genes
most correlated with gene \# 1127 have  greater scores than this gene.
In particular, gene  \#1127 has correlation $>0.6$ with a gene having score
about 4.7. Hence our procedure averages the scores of these two genes
to produce a new score of about 3.8. 

\begin{figure}
\centerline{\epsfig{file=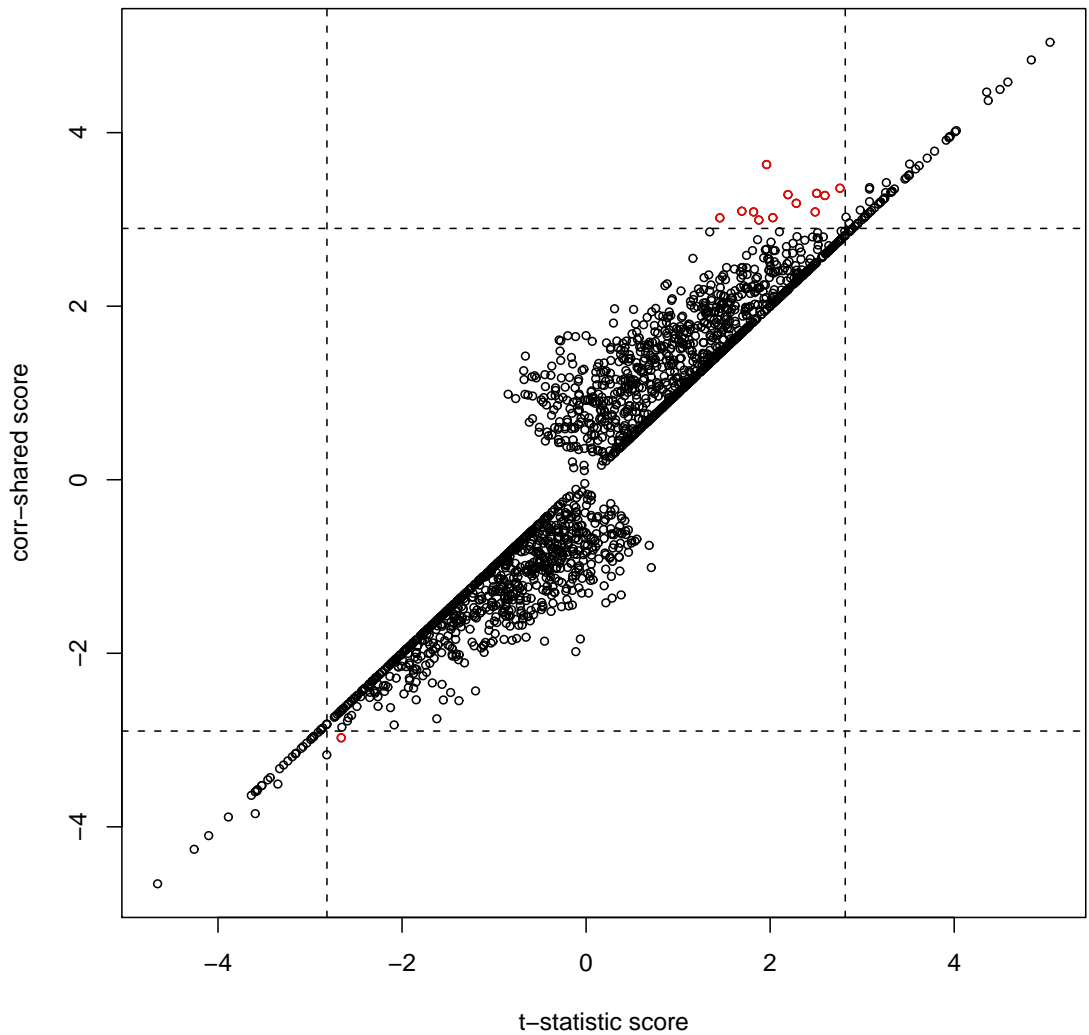,width=6in}}
\caption[skin5]{\small\em Skin data: correlation-shared score versus
t-statistic score. Broken lines are drawn at the cutoffs yielding
100 significant genes for each method. 
The red points are the
the genes that are significant by correlation-sharing but not by 
t-statistic.}
 \label{skin5}
\end{figure}

Figure \ref{skin5}
shows the  correlation-shared score versus the
t-statistic score. Setting the cutoffs so that each method yields
100 significant genes, there are 13 genes which are called  by
each method and not called by the other.
The red points represent
the genes that are called  significant by correlation-sharing but not by the
t-statistic. Many of these genes are highly correlated with each
other, and hence they boost up each other's score.

\begin{figure}
\centerline{\epsfig{file=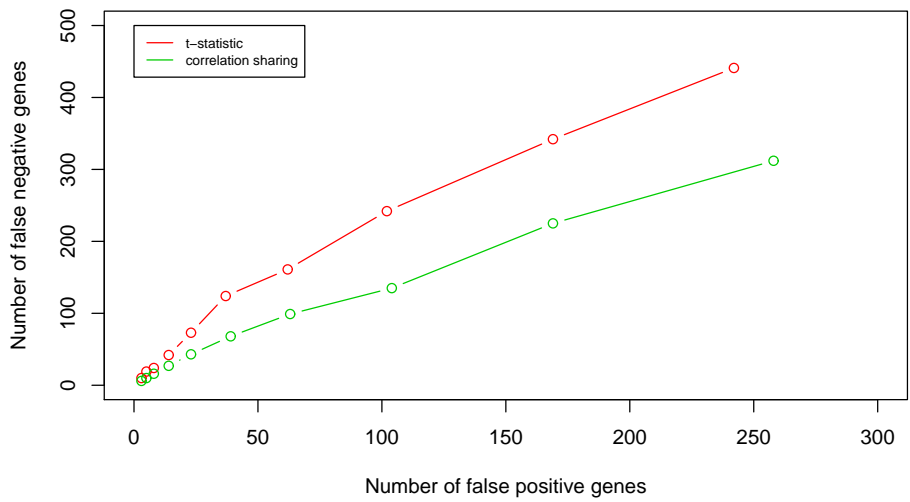,width=5in}}
\caption[skintest]{\small\em Skin data test set results.
Here we formed cutofff rules on the training set, and assessed genes in
in a separate test set.
Shown are the number of false positive and negative genes in the test set,
 as the cutpoint is varied.}
\label{skintest}
\end{figure}

In Figure \ref{skintest} we do another test of our procedure.
We randomly divided the samples into equal-sized training and test sets.
We computed the t-statistic and correlation sharing statistics
on the training set, and also evaluated on the test set.
For each trial cutpoint applied to  the training set scores, we counted the
number of genes with scores above or below this cutpoint 
in the test set. Genes above the cutpoint in the training set but
below it in the test set were considered ``false positives'',
and conversely for false negatives. The results  in Figure
\ref{skintest} show that  correlation sharing has fewer  false negatives
for the same number of false positives.

\section{Example: protein mass spectrometry}
\label{sec.protein}
This example (taken from \cite{carlson2005}) consists of the  intensities of 3160 peaks on 20 patients:
10 healthy patients and 10 with Kawasaki's disease.
They were measured on a SELDI protein mass spectrometer.

\begin{figure}
\begin{psfrags}
\psfrag{# called}{\footnotesize \hskip -.35in Number of genes called significant}
\psfrag{#FP}{\footnotesize \hskip -.35in Number of false positive genes}
\centerline{\epsfig{file=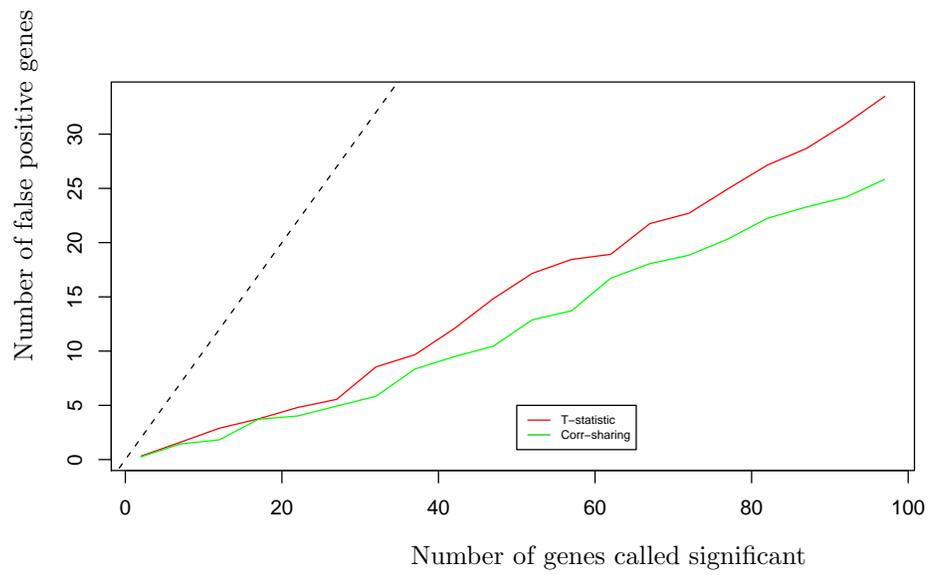,width=5in}}
\end{psfrags}
\caption[cohen1]{\small\em  Results for protein mass spectrometry example}.
\label{cohen1}
\end{figure}

Figure \ref{cohen1} shows that correlation sharing offers a mild improvement in
the false positive rate.

\begin{figure}
\centerline{\epsfig{file=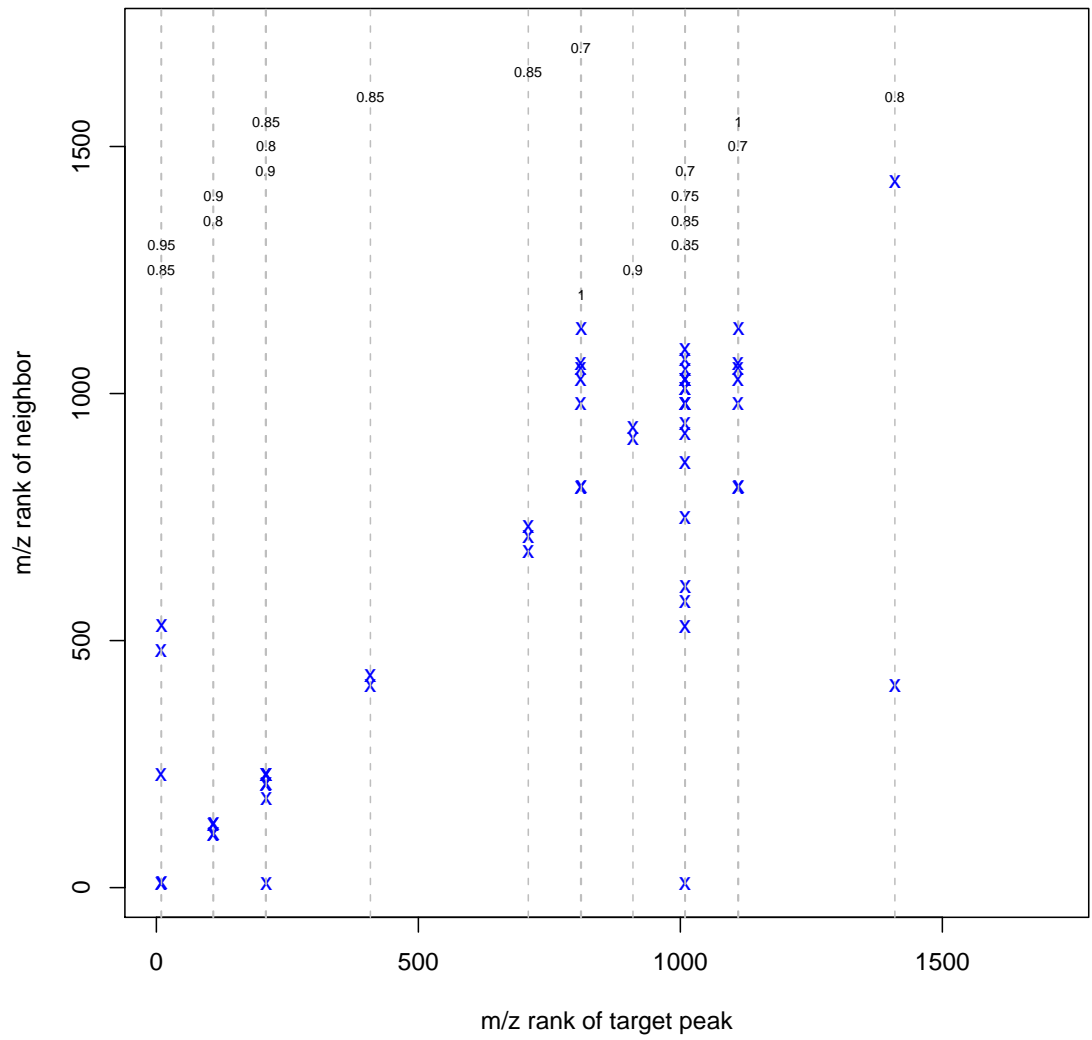,width=6in}}
\caption[cohen2]{\small\em  Protein mass spectrometry example: locations of neighbors of top 50 peaks, for those peaks that were given    given  neighborhoods of more than a single feature. The maximizing correlations are indicated at the
 top of the plot (note that in some cases there are multiple target peaks
 shows near the same position.}
\label{cohen2}
\end{figure}

For the 50 peaks having the top scores,
19 of these peaks were given  neighborhoods of more than a single feature by the correlation
sharing procedure.
The smallest correlation  chosen for neighborhood averaging was 0.7.
Now in this example, each peak has an associated  $m/z$ (mass over charge) location:
this was not used in the correlation-sharing procedure, but we can look posthoc
at the these values within each averaging neighborhood.
Figure \ref{cohen2} shows  the
location of the  each of the 19 peaks (horizontal axis) and the chosen neighbors
(vertical axis). The corresponding neighborhood correlation
is indicated along the top of the plot.
We see that most often, the selected neighbors 
are close to the target peak. But in some cases, they can be very far apart.
Some biological insights might emerge from examination  of these groups of
peaks.

\section{Asymptotic Analysis}
\label{sec.theory}

In this section we show that, under appropriate conditions,
correlation sharing improves power.
More specifically, we show that for null genes,
$U_i$ has similar behavior to $T_i$,
while for nonnull genes,
$U_i$ tends to be stochastically larger than $T_i$.
For simplicity, we focus on a one-sample, one-sided test.
We denote by $X_{ik}$ the measurement for gene $i$ in sample $k$.
Let $T_i=n^{-1}\sum_{k=1}^n X_{ik}$ denote the test statistic for gene $i$ and assume that
$X_{ik}\sim N(\beta_i,\sigma^2)$
where
$\beta_i=0$ for null genes and 
that $\beta_i >0$ for non-nulls.
Let $\rho(i,j)=\corr(X_{ik}, X_{jk})$ denote the true residual correlation between gene $i$ and gene $j$,
$\hat\rho(i,j)$ denote the estimated residual correlation.

The correlation-shared statistic is
\begin{eqnarray}
U_i   &=& \max_{\rho} \frac{1}{|C_\rho(i)|}\sum_{j\in C_\rho(i)} T_j\\
{\rm where}\ C_\rho(i) &=& \Biggl\{j:\ \hat{\rho}(i,j) \geq \rho\Biggr\}.
\end{eqnarray}
Throughout this section we make a small modification to the statistic
which simplifies the analysis:
we restrict the maximization in the definition of $U_i$
to be over correlation neighborhoods no larger than $K$,
where $K$ is some fixed integer.

Recall that there are $m$ genes and $n$ observations.
We require both $m$ and $n$ to grow in the asymptotic analysis.
Typically, $m$ is much larger than $n$
so, to keep the asymptotics realistic, we allow $n$ to grow
very slowly relative to $m$.
Specifically, we assume:
\begin{equation}\label{eq::log-cond}
\hspace{-2cm}{\sf Assumption \ (A1):}\ \ \ \ \ \ \ \ \ \ 
n \equiv n(m) > C \log m\ \ \ {\rm for\ some\ sufficiently\  large\ } C>0.
\end{equation}

Let ${\cal P}$ denote the nonnull genes and let
${\cal N}= {\cal P}^c$ denote the null genes.
We will also need the following:

\vspace{1cm}

\noindent
{\sf Assumption (A2):}
There exist
$0 < \delta <  1$ such that
\begin{equation}\label{eq::delta}
0 = \max_{ \stackrel{i\in {\cal N}}{j\in {\cal P}}} \rho(i,j) =
 \max_{ i,j\in {\cal N}} \rho(i,j) < \delta = \mathop{\min}_{ i,j\in {\cal P}} \rho(i,j).
\end{equation}
Thus we make the strong assumption that there is positive residual correlation among the non-null genes, but no residual correlation among the null genes or between the non and
non-null genes. This simplifies our analysis.
Later, we will relax this assumption.

\vspace{1cm}

{\sf LEMMA 1.}
Assume that (A1) holds.
Fix $\epsilon>0$.
Then, for all large $m$,
\begin{equation}
\max_{ij}|\hat{\rho}(i,j) - \rho(i,j)| < \epsilon\ \ \ a.s.
\end{equation}
and
\begin{equation}
\max_{i}|T_i - \beta_i| < \epsilon\ \ \ a.s.
\end{equation}
That is,
$\hat{\rho}(i,j) = \rho(i,j) + o(1)$, uniformly over $i,j$, a.s.
and
$T_i = \beta + o(1)$, uniformly over $i$, a.s.

\vspace{1cm}

{\sf PROOF of Lemma 1.}
Kalisch and B\"{u}hlmann (2005) show that, 
\begin{equation}
\mathbb{P}(|\hat\rho(i,j)- \rho(i,j)|>\epsilon) \leq 
c_1 (n-1) \exp\left\{ -(n-3)\log ((4+\epsilon^2)/(4-\epsilon^2))\right\}
\end{equation}
for some $c_1>0$.
So,
\begin{eqnarray}
\mathbb{P}(\max_{i,j}|\hat\rho(i,j) - \rho(i,j)| > \epsilon) &&\\
&& \hspace{-5cm}\leq 2 m^2 c_1 (n-1)\times \exp\left\{ -(n-3)\log ((4+\epsilon^2)/(4-\epsilon^2))\right\}
\leq \frac{1}{m^\beta}
\end{eqnarray}
where
\begin{equation}
\beta = \frac{C}{2} \log ((4+\epsilon^2)/(4-\epsilon^2)) -
\frac{\log (4c_1) - \log C - \log\log m}{\log m} -2 >0.
\end{equation}
For $C$ sufficiently large, $\beta>1$.
The first result then follows from the Borel-Cantelli Lemma.
For the second result, apply Mill's inequality:
\begin{equation}
\mathbb{P}(\max_i | T_i - \beta_i | > \epsilon) \leq m e^{-n \epsilon^2/2\sigma^2}.
\end{equation}
The result follows from
assumption (\ref{eq::log-cond}) and the Borel-Cantelli Lemma.
$\blacksquare$

\vspace{1cm}

{\sf LEMMA 2.}
Assume (A1) and (A2).
Then, for all $\rho> \delta$ and each $i\in {\cal P}$,
\begin{equation}
C_\rho(i) \cap {\cal N} = \emptyset\ \ \ \ a.s.
\end{equation}
for all large $m$.
Also, for every $i\in {\cal N}$,
\begin{equation}
C_\rho(i) \cap {\cal P} = \emptyset\ \ \ \ a.s.
\end{equation}
for all $\rho>0$.
Thus, there are no nulls in the correlation neighborhoods
$C_\rho(i)$ of a non-null gene, except possibly for small $\rho$.
Similarly,
there are no nonnulls in the correlation neighborhoods
$C_\rho(i)$ of a null gene.

\subsection{The Oracle Statistic}

To understand the behavior of the correlation sharing statistic, it is helpful
to first consider an oracle version of the statistic
based on the true correlations.
Let
\begin{eqnarray}
\kappa_i    &=& \max_{\rho} \frac{1}{|\nu_\rho(i)|}\sum_{j\in C_\rho(i)} T_j\\
{\rm where}\ \ \ \ \nu_\rho(i) &=& \Biggl\{j:\ \rho(i,j) \geq \rho\Biggr\}.
\end{eqnarray}

Let us fix some nonnull gene $i\in {\cal P}$ and
without loss of generality, take $i=1$.
Without loss of generality, label the genes so that
\begin{equation}
\rho(1,2) >  \rho(1,3) > \cdots > \rho(1,m).
\end{equation}
Then,
\begin{eqnarray}
\kappa_1    &=& \max_r \frac{1}{r} \sum_{i=1}^r T_i 
= \max_r \frac{1}{r} \sum_{i=1}^r \frac{1}{n}\sum_{j=1}^n X_{ij}\\
&=& \max_r \frac{1}{r} \sum_{i=1}^r \frac{1}{n}\sum_{j=1}^n (\beta_i + \epsilon_{ij})
= \max_r \Biggl( \beta(r)  + \frac{1}{rn} \sum_{i=1}^r \sum_{j=1}^n \epsilon_{ij}\Biggr)\\
&=& \max_r \Biggl(\overline{\beta}(r) + Z(r)\Biggr)
\end{eqnarray}
where
\begin{equation}
\overline{\beta}(r) = \frac{1}{r}\sum_{i=1}^r \beta_i
\end{equation}
is the Cesaro average and
$Z(\cdot)$ is a mean zero Gaussian process
with covariance kernel
\begin{equation}\label{eq::cov}
J(r,s) = \frac{\sigma^2}{nrs}\sum_{j=1}^r\sum_{k=1}^s \rho(j,k).
\end{equation}
The distribution of $\kappa_1$ is thus
the distribution of the maximum of a noncentered, nonstationary Gaussian process.

If $\overline{\beta}(r)$ is strongly peaked around some value $r_*$, then
\begin{equation}
\kappa_1  \approx \overline{\beta}(r_*) + Z(r_*) \equiv V_* \sim N(\overline{\beta}(r_*),J(r_*,r_*)).
\end{equation}
Hence,
\begin{equation}\label{eq::bound}
\mathbb{P}(\kappa_1 > t) \approx  \mathbb{P}(V_* > t).
\end{equation}
In particular, suppose that
$\rho(1,i) = \rho$ for $i\in {\cal P}$ and
$\rho(1,i) = 0$ for $i\in {\cal N}$.
Then,
\begin{equation}
V_* \sim N\left( \overline{\beta}(r_*), \frac{1+2\rho}{r_* n}\right)
\end{equation}
and so
\begin{equation}
\mathbb{P}(\kappa_1 > t) \geq  \mathbb{P}\left(\chi_1^2(\pi) > \frac{r_* n t^2}{1+2\rho}\right)
\end{equation}
where
$\chi_1^2(\pi)$ is a noncentral $\chi_1^2$ with noncentrality parameter
\begin{equation}
\pi = \frac{nr_*\overline{\beta}^2(r_*)}{1+2\rho}.
\end{equation}
In contrast, $T_1$ has noncentrality parameter
$n\beta_1^2$.
These heuristics imply that correlation sharing improves the power if
\begin{equation}
\frac{r_* \overline{\beta}^2(r_*)}{1+2\rho} > \beta_1^2,\ \ \ 
{\rm where}\ r_* = {\rm argmax}_r \overline{\beta}(r).
\end{equation}

\begin{figure}
\hspace{2cm}\includegraphics[width=5in]{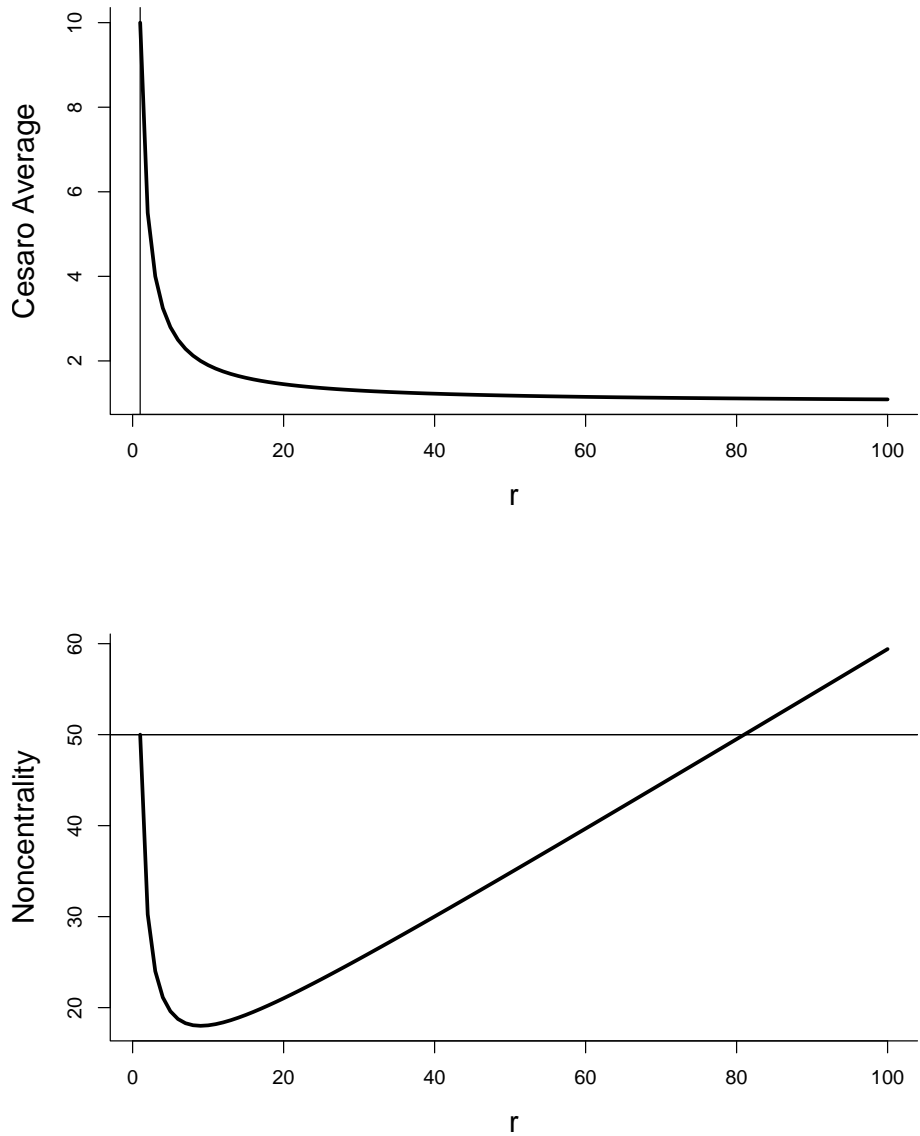}
\caption{\em The non-centrality parameter as a function of neighborhood size; least favorable case.
The top plot is the cumulative average $r^{-1}\sum_{j=1}^r \beta_j$
versus $r$. The bottom plot shows the
noncentrality parameter versus $r$.
The horizontal line shows the noncentrality parameter for $T_1$.
For $1 < r \leq 80$, the noncentrality parameter for $T_1$
is larger than noncentrality parameter for $U_1$.
Since the top plot is maximized at $r=1$ we expect that the correlation neighborhood
for $U_1$ shoule have $r$ close to 1.}
\label{fig::noncen1}
\end{figure}

\begin{figure}
\hspace{2cm}\includegraphics[width=5in]{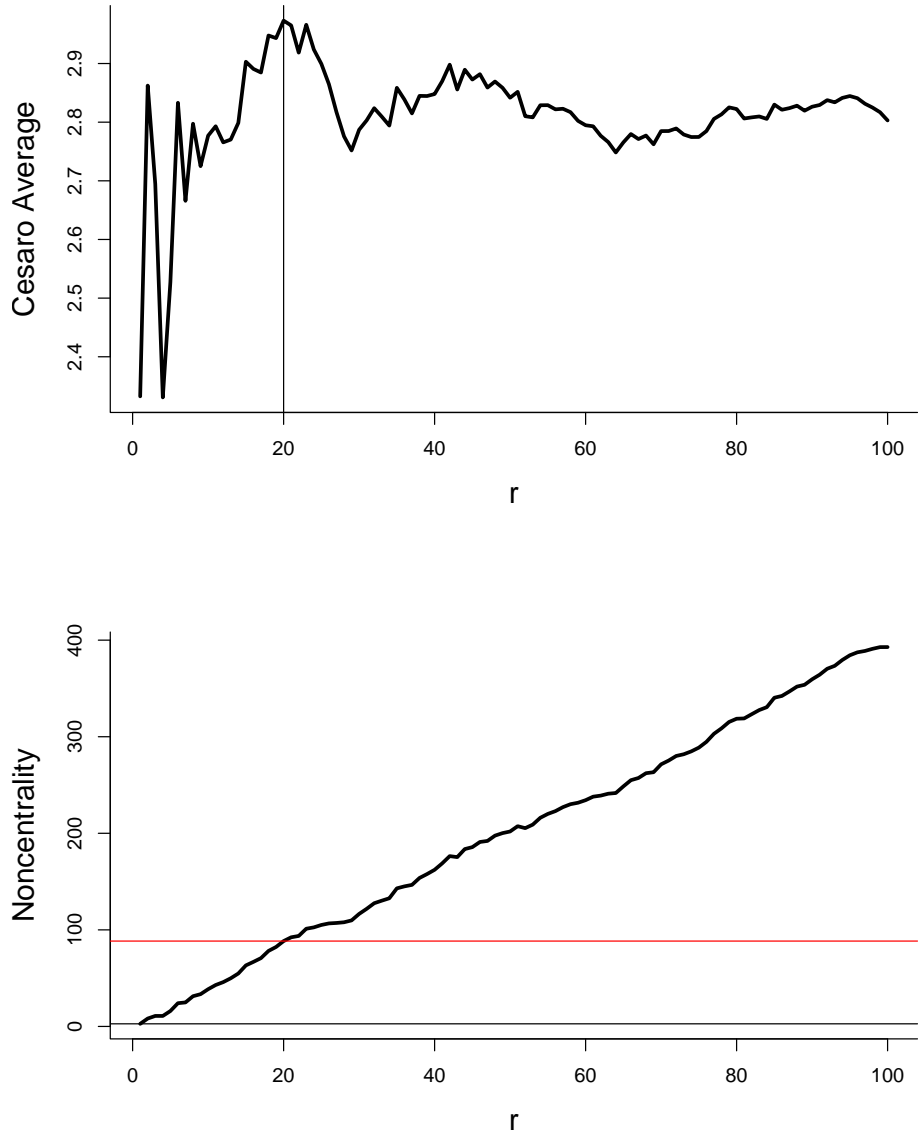}
\caption{\em The non-centrality parameter as a function of neighborhood size; typical case.
The top plot is the cumulative average $r^{-1}\sum_{j=1}^r \beta_j$
versus $r$. The bottom plot shows the
noncentrality parameter versus $r$.
The horizontal line near 0 shows the noncentrality parameter for $T_1$.
The horizontal line near 100 shows the noncentrality parameter for $U_1$
when the correaltion neighborood is $r=20$ corresponding to the maximum of the top plot.
Not only is there a large gain in noncentrality, but the gain is robust to fluctuations in $r$.}
\label{fig::noncen2}
\end{figure}

Figures \ref{fig::noncen1} and \ref{fig::noncen2} 
illustrate this analysis.
The top plot in each figure
is $\overline{\beta}(r)$ and
the bottom plot is the noncentrality as a function of the size $r$ of the
correlation neighborhood.

Figures \ref{fig::noncen1} shows a least favorable case
in which $\beta_1=10$ and $\beta_i =1$ for
$i>1, i\in {\cal P}$.
(In all cases we took $\rho=.5$).
We call this least favorable since 
$T_1$ has the largest mean; any averaging can only reduce
its mean.
Now, $r_* =1$ and $T_1$ has noncentrality 50.
The randomness of $\hat\rho$ can lead to a correlation
neighborhood larger than $r_*=1$.
If so, the noncentrality parameter
can be reduced as is evident from the steep
decline of the curve in the second plot.

Figure \ref{fig::noncen2} shows a more realistic case.
Here we used a random effects model
and took $\beta_i \sim N(3,1)$.
This makes $\sum_r \beta_r$
a random walk.
Correspondingly,
$\overline{\beta}(r)$ behaves like a random walk
for small $r$ but settles down to a constant for large $r$.
In this case, $r_*$ tends to be small but
the noncentrality grows rapidly.
The result is a dramatic gain in noncentrality.
Also, the gain is robust to the choice of $r$.

Now consider a null gene $i\in {\cal N}$.
Again take $i=1$.
Then, by assumption (A2), $\rho(1,j)=0$ for all $j>1$.
Hence, $\nu_\rho(1)=\{1\}$ for all $\rho>0$ and
$\kappa_1    = T_1$
so the null distribution is unaffected by
correlation sharing.

Let us now consider weakening (A2).
Suppose we allow some small, nonzero correlation
$\Delta$ among null genes.
Change the definition of $U_i$ to
\begin{equation}
U_i   = \max_{  \stackrel{\rho>\Delta}{|C_\rho(i)|\leq K}  } 
\frac{1}{|C_\rho(i)|}\sum_{j\in C_\rho(i)} T_j
\end{equation}
Now replace (A2) with:

\vspace{1cm}

\noindent
{\sf Assumption (A2'):}
\begin{equation}
\mathop{\min}_{ \stackrel{i\in {\cal P}}{j\in {\cal P}}} \rho(i,j) >
\max_{ \stackrel{i\in {\cal N}}{j\in {\cal P}}} \rho(i,j)
\end{equation}
and 
\begin{equation}
\max_{ \stackrel{i\in {\cal N}}{j\in {\cal P}}} \rho(i,j) < \Delta.
\end{equation}

The analysis for nonnull genes is virtually unchanged.
For null genes, 
condition (A2') ensures that
$\kappa_1 = T_1$.
An interesting extension is to estimate $\Delta$
from the data.
We leave this to future work.

%

\subsection{Relationship Between $U_i$ and the Oracle}

The analysis in the previous section ignores the
variability of the $\hat{\rho}(i,j)'s$.
Now we relate $\kappa_i$ to $U_i$.

First, under appropriate assumptions, we will show that
for nonnull genes,
$U_1$ is at least as large as $\kappa_1$.
Suppose there exists a
decreasing function $f:[0,1] \to [0,1]$ with
$f(0)=1$, such that
\begin{equation}
\rho(1,i) = f(i/m).
\end{equation}

Suppose that $f$ is
a simple function, that is, $f$ takes
finitely many values
$a_1 > a_2 > \cdots > a_k$.
The level sets
$\nu_\rho(1) = \{ j:\ \rho(1,j) \geq \rho\}$
can only be of the form
${\cal A}_s = \{ j:\ \rho(1,j) \geq a_s\}$
for $s=1, \ldots, k$.
Choose $\epsilon >0$ small.
By Lemma 1, $\max_j |\hat{\rho}(1,j) - \rho(1,j)| < \epsilon$ a.s.
Let $I= \{ \rho\in[0,1]:\ \min_s |\rho - a_s| > \epsilon\}$.
For all $\rho\in I$,
$\nu_\rho(1) = C_\rho(1)$ a.s.
Then, for all large $m$,
\begin{eqnarray}
U_1   &=& \max_{\rho} \frac{1}{|C_\rho(1)|}\sum_{j\in C_\rho(1)} T_j
\geq  \max_{\rho\in I} \frac{1}{|C_\rho(1)|}\sum_{j\in C_\rho(1)} T_j\\
&=& \max_{\rho\in I} \frac{1}{|\nu_\rho(1)|}\sum_{j\in C_\rho(1)} T_j\ \ a.s.
= \max_{\rho} \frac{1}{|\nu_\rho(1)|}\sum_{j\in C_\rho(1)} T_j\\
&=& \kappa_1
\end{eqnarray}
so that $U_1$ is at least as large as $\kappa_1$.

Now we drop the assumption that
$f$ is simple and instead assume it is
continuous and strictly decreasing.
Similarly, assume there
exists a continuous, integrable
function $g$
such that
\begin{equation}
\beta_i  = g(i/m).
\end{equation}
Suppose that
$\overline{g}(u)\equiv s^{-1}\int_0^s g(u) du$ is maximized at some $s_* >0$.
Let $c = f(s_*)$ and $r = |\nu_c(1)|$.
Then, a.s. for all large $m$,
\begin{eqnarray}
\kappa_1    &=& \max_{\rho} \frac{1}{|\nu_\rho(1)|}\sum_{j\in \nu_c(1)} T_j=
\max_{\rho} \frac{1}{|\nu_\rho(1)|}\sum_{j\in \nu_\rho(1)} \beta_j + o(1)\\
&=& \max_r \frac{1}{r} \sum_{j=1}^r \beta_j + o(1)
= \max_s \frac{1}{s} \int_0^s g(u) du + o(1)\\
&=& \frac{1}{s_*} \int_0^{s_*} g(u) du + o(1) = \overline{g}(s_*) + o(1).
\end{eqnarray}
Hence,
\begin{equation}
\kappa_1 \approx  \frac{1}{|\nu_c(1)|}\sum_{j\in \nu_\rho(i)} T_j=
 \frac{1}{r} \sum_{j:\ \rho(1,j)\geq c}T_j .
\end{equation}
Let $R = |C_c(1)|$.
From Lemma 1 and the assumptions on $f$,
$R/m=r/m +o(1)$ a.s. and
\begin{eqnarray}
U_1   &=& \max_{\rho} \frac{1}{|C_\rho(1)|}\sum_{j\in C_\rho(1)} T_j
 \geq  \frac{1}{|C_c(1)|}\sum_{j\in C_c(1)} T_j\\
&=& \frac{1}{R}\sum_{j:\ \hat\rho(1,j)\geq c} T_j\\
&=& \frac{1}{R}\sum_{\rho(1,i)\geq c} T_i +
\frac{1}{R}\sum_{\stackrel{\hat\rho(1,i)\geq c}{\rho(1,i)< c}} T_i -
\frac{1}{R}\sum_{\stackrel{\hat\rho(1,i)< c}{\rho(1,i)\geq c}} T_i \\
&=& \frac{1}{r}\sum_{\rho(1,i)\geq c} T_i + o(1)\\
&=& \kappa_1.
\end{eqnarray}
Thus,
$U_1 \geq \kappa_1 + o(1)$.

Now suppose that $i=1$ is a null gene.
Fix a small $\epsilon >0$.
Under (A2), we eventually, have
\begin{equation}
| \{j>1 :\ \hat\rho(1,j) > \epsilon \}| =0
\end{equation}
and hence
\begin{equation}
U_1 = \kappa_1 \ \ a.s.
\end{equation}
The same holds under (A2').

\section{Other issues}
\label{sec.other}
Computation of the correlation shared statistic can be
challenging when the number of features $m$ is large.
Brute force computation is $O(m^2)$. 
In principle, a  KD tree   can be used to quickly find the neighbors of a given point
with correlation at least $\rho$. The building of the tree requires $O(m\log m)$
computations, while the nearest neighbor search takes $O(\log{m})$ computations.
Hence the nearest neighbor search for all points requires $O(m\log m)$
computations.
However, since the dimension of the feature space ($n$) is large
n these problems (at least 50 or 100), the KD tree approach is not likely
to be effective in practice (J. Friedman, personal communication).

Hence we instead do a direct brute force computation, exploiting the sparsity
of the set of pairs of points with large correlation. The resulting procedure
is quite fast, requiring for example 2.7s on the proteomics example
($m=3160, n=20$).

The proposal of this paper can be applied to outcome
measures  other than two-class problems.
We have seen  this earlier in the lymphoma example, where the outcome
was survival time.
Other response types that may arise include a  multi-class or quantitative
outcome.
The modification to the correlation-sharing  technique is simple: the t-statistic (\ref{tstat}) is simply replaced by
a score
that is appropriate for the outcome measure. 
For survival data, for example, we use the partial likelihood score statistic for
each gene. This was illustrated in the lymphoma data of
Table \ref{tab:summary}.

Correlation-sharing provides a recipe for supervised clustering of features.
Hence one might use correlation-sharing as a pre-processing step,
by averaging the given features in the prescribed clusters.
Then these averaged features could be used as input into a regression
or classification procedure. This is a topic for future study.

\bigskip

{\bf Acknowledgments} 
We would like to thank John Storey for 
showing us a pre-preprint of his ``optimal discovery procedure'' paper.
We would also  like to thank Jerry Friedman for helpful discussions.
Tibshirani was partially supported by National Science Foundation
Grant DMS-9971405  and National Institutes of Health Contract  N01-HV-28183.

\bibliographystyle{agsm}

\bibliography{/home/tibs/texlib/tibs}

\begin{table}
{\begin{small}
\begin{tabular}{lllll}
Name &  Description&  \# Samples &\# Features & Source \\
\hline
Skin &  Two classes& 58 & 12,625  & \citeasnoun{jK01}\\
Duke breast cancer &  Two classes & 49 &  7097  & \citeasnoun{huang2003}\\
BRCA &  Two classes& 15 & 3226  & \citeasnoun{Hedetal2001}\\
Lymphoma &     Survival& 240 & 7399  & \citeasnoun{aR02}\\
\end{tabular}
\caption[tab:summary]{\em Summary of datasets for Figure \ref{cancer.corr}.}
\label{tab:summary}
\end{small} }
\end{table}

\begin{figure}
\centerline{\epsfig{file=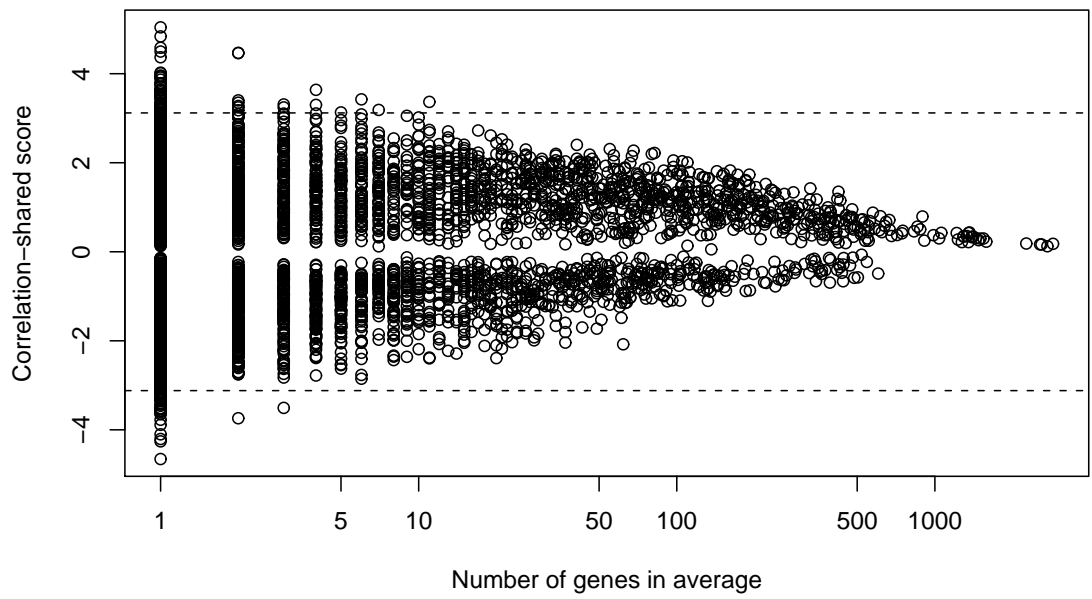,width=6in}}
\caption[skin4]{\small\em Skin data: correlation-shared score versus 
number of genes used in each gene average; horizontal lines are drawn at cutpoints
that yield  100 significant genes. Note that most of the  significant genes 
use no averaging, and none use a window of more than 10 genes} 
 \label{skin4}
\end{figure}

\end{document}